\documentclass{article}
\usepackage{amsmath, amssymb, latexsym}
\title{More Constructions for Tur\'{a}n's (3, 4)-Conjecture}
\author{Andrew Frohmader}

\newtheorem{theorem}{Theorem}[section]

\newtheorem{lemma}[theorem]{Lemma}
\newtheorem{definition}[theorem]{Definition}
\newtheorem{conjecture}[theorem]{Conjecture}

\newtheorem{construction}[theorem]{Construction}
\newtheorem{problem}[theorem]{Problem}

\def\proof{\smallskip\noindent {\it Proof: \ }}
\def\endproof{\hfill\ensuremath{\square}\medskip}

\begin{document}

\maketitle

\begin{abstract}
For Tur\'{a}n's (3, 4)-conjecture, in the case of $n = 3k+1$ vertices, ${1 \over 2}6^{k-1}$ non-isomorphic complexes are constructed that attain the conjecture.  In the case of $n = 3k+2$ vertices, $6^{k-1}$ non-isomorphic complexes are constructed that attain the conjecture.
\end{abstract}

\section{Introduction}

In 1941, Tur\'{a}n \cite{turan} posed a problem about edges of hypergraphs.  Translated into the language of simplicial complexes, it reads as follows.  Suppose that a simplicial complex with $n$ vertices has every minimal non-face on exactly $k$ vertices, and has no faces on $r$ vertices.  How many faces on $k$ vertices can it have, as a function of $n$, $k$, and $r$?

Let $t_k(n, r)$ be the greatest number of faces on $k$ vertices that the complex can have.  If all possible sets of $k$ vertices formed a face, there would be ${n \choose k}$ faces on $k$ vertices.  Hence, $t_k(n, r) / {n \choose k}$ is the fraction of the potential faces on $k$ vertices.  If a complex on $n+1$ vertices attains $t_k(n+1, r)$, then removing one vertex leaves a complex with $n$ vertices and at most $t_k(n, r)$ faces on $k$ vertices.  Average over all the ways to remove a vertex and we get $${t_k(n+1, r) \over {n+1 \choose k}} \leq {t_k(n, r) \over {n \choose k}}.$$

Tur\'{a}n then posed the following problem.

\begin{problem}  \label{genproblem}
Fix $r > k$.  Let $C$ be a simplicial complex on $n$ vertices with no face on $r$ vertices.  Suppose further that every minimal non-face of $C$ has $k$ vertices.  Let $t_r(n, k)$ denote the greatest number of faces on $k$ vertices that the complex $C$ could possibly have. Compute
$$\lim_{n \to \infty} {t_k(n, r) \over {n \choose k}}.$$
\end{problem}

As we have seen, the limit is of a (weakly) decreasing sequence of positive numbers, so it must converge.  Tur\'{a}n's theorem \cite{turan} established that if $k = 2$, the answer is ${r-1 \over r}$, but this is the only case where the answer is known.

The next simplest case is when $k = 3$ and $r = 4$.  For this, Tur\'{a}n \cite{turanconj} conjectured the following.

\begin{conjecture}  \label{turanconj}
Let $C$ be a simplicial complex of dimension two, for which every minimal non-face has three vertices.  Let the number of vertices of $C$ be $n$.  Then the number of triangles of $C$ is at most
$$\left\{
\begin{array}{ccl}
{5 \over 2} k^3 - {3 \over 2} k^2 & \textup{if} & n = 3k \\ {5 \over 2} k^3 + k^2 - {1 \over 2} k & \textup{if} & n = 3k+1 \\ {5 \over 2} k^3 + {7 \over 2} k^2 + k & \textup{if} & n = 3k+2.
\end{array}
\right.$$
\end{conjecture}

Conjecture~\ref{turanconj} would imply an answer of ${5 \over 9}$ to this case of Problem~\ref{genproblem}.  Tur\'{a}n established ${5 \over 9}$ as a lower bound by giving the following construction that attains the bound of his conjecture.  Divide the $n$ vertices into three parts as evenly as possible, and arrange the parts cyclically so that each has one to its ``left" and one to its ``right".  The complex has all $n$ vertices and all ${n \choose 2}$ possible edges.  The triangles of the complex are those for which one vertex is from each part, or two are from one part and one from the part to its right.

For $n \geq 7$, this is not the only construction that attains the conjecture, however.  Brown \cite{brown} showed that there are at least $k-1$ non-isomorphic constructions that attain the bound if $n = 3k$.  Kostochka \cite{kostochka} generalized Brown's constructions to give $2^{k-2}$ non-isomorphic constructions if $n = 3k$.  These constructions are easiest to describe in terms of which triangles are not in the complex, and Conjecture~\ref{turanconj} can be reformulated as a lower bound on the number of missing triangles, given by ${n \choose 3}$ minus the formulas in the conjecture.

Kostochka further observed that by removing one or two vertices from his constructions, one could obtain many constructions that attain the bound of Conjecture~\ref{turanconj} if $n$ is not a multiple of 3.  Removing $j$ vertices can give on the order of $k^j2^k$ constructions, but many of them are isomorphic to each other or do not attain the bound.  This paper improves on that result by showing that there are on the order of $6^k$ non-isomorphic constructions that attain the bound of Conjecture~\ref{turanconj} if $n = 3k+1$ or $n = 3k+2$.

\begin{theorem} \label{mytheorem}
If $n = 3k+1$, then there are at least ${1 \over 2}(6)^{k-1}$ complexes that attain the bound of Conjecture~\ref{turanconj}, no two of which are isomorphic.  If $n = 3k+2$, then there are at least $6^{k-1}$ complexes that attain the bound of Conjecture~\ref{turanconj}, no two of which are isomorphic.
\end{theorem}

Some upper bounds for the $k = 3, r = 4$ case of Problem~\ref{genproblem} are also known.  In particular, if we can compute $t_3(n, 4)$ for any particular value of $n$, then we get $t_3(n, 4) / {n \choose k}$ as an upper bound on the limit. Some better upper bounds were given by de Caen \cite{decaen0} of $\approx .6213$, Giraud (unpublished, see \cite{decaen}) of $(\sqrt{21} - 1) / 6 \approx .5971$, and Chung and Lu \cite{chunglu} of $(3 + \sqrt{17}) / 12 \approx .5936$. Conjecture~\ref{turanconj} has been verified for the cases $n \leq 13$ by Spencer \cite{spencer}.

The layout of this paper is as follows.  In Section~2, we give our construction.  We start with a more general construction that has some nice properties but is too broad, and then restrict it to the complexes we actually want.  The remainder of the section shows that the complexes of Construction~\ref{subconst} really do attain the bound of Conjecture~\ref{turanconj}, and counts the number of such complexes.  Section~3 shows that no two of the complexes of Construction~\ref{subconst} are isomorphic to each other.  In Section~4 we discuss whether there are complexes other than those of Construction~\ref{subconst} that attain the bound of Conjecture~\ref{turanconj}.

\section{The construction}

This section starts by giving a general construction for which it is clear that the complexes satisfy the requirements of Problem~\ref{genproblem}.  This construction includes many complexes that do not attain the bound of Conjecture~\ref{turanconj}, as well as many complexes that are isomorphic to each other.  We then narrow this repeatedly, culminating in Construction~\ref{subconst}, a subset of the complexes of the first construction.  This last construction eliminates the complexes that do not attain the bound of the conjecture and all but one complex from each isomorphism class.  This section also shows that all of the complexes of Construction~\ref{subconst} do attain the bound of Conjecture~\ref{turanconj} and counts how many such complexes there are.

When narrowing down the construction, we sometimes give only the main ideas of the proof rather than a complete proof.  This is because we are mainly interested in showing that there are many complexes that attain the bound of the conjecture without being isomorphic to each other, not that many complexes are isomorphic to each other or do not attain the bound.

\begin{construction} \label{mainconst}
\textup{Divide $n$ vertices into three columns.  Arrange the columns cyclically so that each has one to its ``right" and one to its ``left"; if you start at one column and go to its right three times, you end up back at the original column.  Put a total order on the vertices.  We can refer to one vertex being higher or lower than another in this order.  Color each vertex either red or blue.}

\textup{Construct a simplicial complex with the $n$ vertices as its vertex set, all ${n \choose 2}$ possible edges, and no faces on four or more vertices.  The triangles in the complex are all possible sets of three vertices except for}
\begin{enumerate}
\item \textup{three vertices in the same column, with the top two the same color;}
\item \textup{two vertices in one column, with the higher of the two red, and one vertex in the column to its right;}
\item \textup{two vertices in one column, with the higher of the two blue, and one vertex in the column to its left; and}
\item \textup{two vertices in one column and one vertex in a different column, with the highest vertex in the left column blue, the highest in the right column red, and the two lowest vertices in the same column.}
\end{enumerate}
\end{construction}

We can check cases to show that, for any four possible vertices, some three of them do not form a triangle.  This construction does not always attain the bound of Tur\'{a}n's (3, 4)-conjecture, however.

Observe that the relative order of vertices in different columns matters only if we are considering missing triangles of type 4.  If we have a blue vertex in one column and a red vertex in the column to its right, these two vertices form a missing triangle with each vertex in the column with the lower of the two and below that vertex.  To maximize the number of triangles, the lower of the first two vertices should be whichever has fewer other vertices below it in its column.  We can thus arrange the vertices into rows as well as columns, with a vertex in a higher row always higher than a vertex in a lower row.  If arranging the vertices into rows cannot give a new complex isomorphic to the previous one, then the new one will have more triangles.

We can also require the vertices to be distributed among the columns as evenly as possible, as any construction which does not do this either is isomorphic to one that does or else does not attain the bound.  An argument analogous to that of Theorem~\ref{optcolor} can show that if there were any exceptions to this, there would have to be an exception which leaves a column completely empty.

If there is such a new complex that matches the bound, then there must be one that is not isomorphic to any construction with the vertices distributed more evenly among the columns.  If the top vertex is red, we can shift it to the right one column and change it to blue and get a complex isomorphic to what we began with.  The reverse of this can also be done if the top vertex is blue.  After arranging the vertices into rows, if this isomorphism makes the vertices distributed more evenly, it contradicts our choice of the complex.

If this isomorphism does not make the vertices more evenly distributed, then unless all vertices in the column with the top vertex are the same color, one can show that swapping the column's highest vertex of the opposite color from the top vertex with the one above it increases the number of triangles in the complex.  If all vertices in the column are the same color, then at most half of the possible sets of three vertices form a triangle, so the complex does not reach the bound of Conjecture~\ref{turanconj}.

Adding these restrictions to Construction~\ref{mainconst} gives us the following.

\begin{construction} \label{somecond}
\textup{Divide $n$ vertices into 3 columns and $\big\lceil {n \over 3} \big\rceil$ rows, such that each choice of a column and row has at most one vertex.  All empty spots must be in the top row.  Put a total order on the vertices such that a vertex in a higher row is always greater than a vertex in a lower row.  Color each vertex either red or blue.  The faces of the complex are precisely as given in Construction~\ref{mainconst}.}
\end{construction}

Some of these constructions still do not attain the bound of Conjecture~\ref{turanconj}.  We place some restrictions on colorings to get the following construction as a subset of the complexes of Construction~\ref{somecond}.

\begin{construction} \label{morecond}
\textup{Among the complexes of Construction~\ref{somecond}, eliminate all except those which satisfy the following coloring conditions.  If there are $n = 3k$ vertices, then each row must have all three of its vertices the same color.}

\textup{If there are $n = 3k+1$ vertices, then for all choices of $j \leq k$, if we restrict to the top $j$ rows, the number of red vertices in a column is not more than in the column to its left, except that the column with the top vertex may have one more red vertex than the column to its left.}

\textup{If there are $n = 3k+2$ vertices, for all choices of $j \leq k$, if we restrict to the top $j$ rows, the number of red vertices in a column is not fewer than in the column to its left, except that the column without a vertex in the top row may have one red vertex fewer than the column to its left.}
\end{construction}

The intuitive idea of the coloring conditions is that the red vertices and the blue vertices must each be distributed among the columns as evenly as possible throughout the complex.

An equivalent explanation of the coloring condition if $n = 3k+1$ is that you can hit all of the red vertices by starting at the highest red vertex in the column with the top vertex, and jumping to the next each time by moving one column to the right and possibly down some number of rows, but not up.  If $n = 3k+2$, you move left one column each time instead of right, and start in the right column of the two with a vertex in the top row.

\begin{definition}
\textup{A \textit{color set} of vertices in Construction~\ref{mainconst}, \ref{somecond}, \ref{morecond}, or \ref{subconst} consists of all red vertices from one column and all blue vertices from the column to its right.  The \textit{size} of a color set is its number of vertices.  We say that a vertex is \textit{below a color set} if it is lower than all vertices of the color set in the same column, even if it is higher a vertex of the color set in the other column.}
\end{definition}

Note that all of the colored vertices in a complex are partitioned into three color sets.  The next lemma states that if we remove some rows from the bottom of Construction~\ref{morecond}, the remaining vertices are divided into color sets as evenly as possible.

\begin{lemma}  \label{redbluesplit}
Let $k = \big\lfloor{n \over 3}\big\rfloor$.  If the bottom $j$ rows of a complex from Construction~\ref{morecond} are removed for some $1 \leq j \leq k$, then the total number of vertices in each color set is either $k-j$ or $k-j+1$.
\end{lemma}

\proof  If $n = 3k$, then a color set has exactly one vertex in each row.  As such, if we remove the bottom $j$ rows, each such set has $k-j$ vertices.

Otherwise, let the number of red vertices be $3p+q$ with $0 \leq q \leq 2$.  Let column $A$ be the one with the top vertex if $n = 3k+1$ and the left column with a vertex in the top row if $n = 3k+2$, column $B$ be the column to the right of $A$, and column $C$ be the column to the right of $B$.  One can count the number of red and blue vertices in each column to get the following table.
$$\begin{array}{cccccccc}
\textup{vertices} & q & A \textup{ red} & B \textup{ blue} & B \textup{ red} & C \textup{ blue} & C \textup{ red} & A \textup{ blue} \\ 3k+1 & 0 & p & k-j-p & p & k-j-p & p & k-j-p+1 \\ 3k+1 & 1 & p+1 & k-j-p & p & k-j-p & p & k-j-p \\ 3k+1 & 2 & p+1 & k-j-p-1 & p+1 & k-j-p & p & k-j-p \\ 3k+2 & 0 & p & k-j-p+1 & p & k-j-p & p & k-j-p+1 \\ 3k+2 & 1 & p & k-j-p & p+1 & k-j-p & p & k-j-p+1 \\ 3k+2 & 2 & p+1 & k-j-p & p+1 & k-j-p & p & k-j-p
\end{array}$$

Note that in all cases, the size of each color set is either $k-j$ or $k-j+1$.  \endproof

A converse of this lemma is true, as well:  if a complex of Construction~\ref{somecond} satisfies the color set conditions of Lemma~\ref{redbluesplit}, then it must satisfy the coloring conditions of Construction~\ref{morecond}.  This is readily proven by counting the number of constructions of each type, observing that they are the same, and that one set is a subset of the other by Lemma~\ref{redbluesplit}.

Together with Lemma~\ref{redbluesplit}, the next theorem effectively states that all of the complexes of Conjecture~\ref{morecond} attain the bound of Conjecture~\ref{turanconj}.

\begin{theorem} \label{optcolor}
Let $k = \lfloor{n \over 3}\rfloor$.  A simplicial complex given by Construction~\ref{somecond} attains the bound of Conjecture~\ref{turanconj} if and only if it satisfies the following condition:  for all $j$ with $1 \leq j \leq k$, if the bottom $j$ rows are removed, then each color set has either $k-j$ or $k-j+1$ vertices.
\end{theorem}

\proof  Suppose that apart from the bottom $j$ rows, one column has $a$ red vertices and $b$ blue vertices, the column to its right has $c$ red vertices and $d$ blue, and the column to its right has $e$ red vertices and $f$ blue.  We wish to count the number of missing triangles whose lowest vertex is in the $(j-1)$-th row from the bottom.  For missing triangles of the first type, we get ${a \choose 2} + {b \choose 2} + {c \choose 2} +{d \choose 2} + {e \choose 2} + {f \choose 2}$.  The number of missing triangles of the second type for which the lowest row used has a vertex in the column with two vertices is $a(c+d+1) + c(e+f+1) + e(a+b+1)$.  The number of analogous missing triangles of the third type is $b(e+f+1) + d(a+b+1) + f(c+d+1)$.  The number of missing triangles of the second or third type for which the two vertices in the same column are both in higher rows than the third vertex is ${a+b \choose 2} + {c+d \choose 2} + {e+f \choose 2}$.  The number of missing triangles of the fourth type is $bc + de + af$.

If we add all of these up and rearrange, we get $${1 \over 2}\big((a+d)^2 + (b+e)^2 + (c+f)^2 + (a+b+c+d+e+f)^2\big).$$  Note that the last term does not depend on the coloring, so we only need to minimize the rest of the sum.  The way to minimize the sum of the squares of three integers with a fixed sum ($n - 3j$) is to make the three integers as close to each other as possible.  In other words, each color set must have size either $k-j$ or $k-j+1$.

In order to minimize the number of missing triangles throughout the complex, this must happen for every choice of $j$, so the complexes satisfying the condition given in the lemma are precisely the optimal complexes among those of Construction~\ref{somecond}.  If we color all vertices red to get Tur\'{a}n's original construction, we satisfy the condition and attain the bound of his conjecture.  \endproof

Some possible ways of coloring and ordering vertices still give complexes isomorphic to each other.  We have already discussed the isomorphism from shifting the top vertex by one column and inverting its color.  Another such isomorphism comes from inverting the color of all vertices and reversing the order of the columns.  More come from changing the color of vertices in the bottom row, which are not used to determine whether any triangles exist.  Still more isomorphisms come from swapping the relative order of two vertices in the same row if it does not affect any missing triangles of type 4.  We can place additional conditions on the colorings and impose an equivalence relation to eliminate all but one complex from any isomorphism class given by these isomorphisms.

\begin{construction} \label{subconst}
\textup{Among simplicial complexes of Construction~\ref{morecond}, discard those which have a blue vertex in the top row.  If $n = 3k+1$, discard complexes for which the top vertex in a column is blue, even if not in the top row.  Consider two complexes with the same coloring to be the same complex if they agree on the relative order of pairs of vertices within a row for which the left vertex is blue and the right vertex is red, even if they disagree on the relative order of other pairs of vertices within a row.  Additionally consider two complexes the same if they disagree only in their bottom row.}
\end{construction}

We typically regard vertices in the bottom row as being uncolored.  By abuse of language, however, we sometimes say that the bottom $j$ vertices of a column are red when we really mean, except for the bottom vertex, which is uncolored.

If $n = 3k$, this gives precisely Kostochka's construction from \cite{kostochka}.  Otherwise, we wish to count how many non-isomorphic complexes this construction gives.  There is only one way to color the top (partially filled) row, and essentially only one way to color the bottom row, as its colors do not matter.

For most intermediate rows, we can pick zero, one, two, or three red vertices.  If we pick one or two, at some point, there will be a blue vertex immediately to the left of a red one in the same row.  Which of these two vertices is higher will matter as some complexes of type 4 will use these pairs of vertices.  There are two ways to pick which of the vertices is higher for each of one or two red vertices in a row.  As such, there are six ways to fill in a row.

The exception to this is the second row in the $n = 3k+1$ case. The column to the left of the top vertex is assumed to have a red vertex, and this forces the column to the left of it to have a red vertex.  The remaining vertex can be red or blue, and if blue, higher or lower than the red vertex to its right. This creates three ways to fill in this row.  As such, we have ${1 \over 2}(6)^{k-1}$ complexes if $n = 3k+1$ and $6^{k-1}$ complexes if $n = 3k+2$.

\section{Checking for isomorphisms}

In the previous section, we constructed many complexes that attain the bound of Conjecture~\ref{turanconj}.  A priori, these could include many ways of presenting the same complex, so that there would be far fewer distinct complexes than given. Indeed, Constructions~\ref{mainconst}, \ref{somecond}, and \ref{morecond} did exactly that.  In this section, we show that no two complexes of Construction~\ref{subconst} are isomorphic.  References to any construction in this section are assumed to come from this last construction.

The main idea of the proof is that we start by defining some combinatorial invariants that are clearly preserved by isomorphisms.  If some invariant differs for two complexes, then they cannot be isomorphic to each other.  We then show how to easily compute the invariants in most cases from the structure of Construction~\ref{subconst}.  This leads to a way to essentially pick out the bottom vertex of each column as given in the construction.  We can then show that no two complexes are isomorphic by removing the bottom row of a complex and proceeding by induction.

Our proof is similar in flavor to Kostochka's proof that no two of his complexes are isomorphic in \cite{kostochka}.  Our proof is somewhat more complicated, as Kostochka was able to exploit all three columns being identical, which is untrue of Construction~\ref{subconst}.  Throughout this section, we assume that the complexes have $n = 3k+1$ or $n = 3k+2$ vertices, though the lemmas would generally hold (sometimes with slight modifications) for $n = 3k$ vertices. First, we need the combinatorial invariants.

\begin{definition}
\textup{An \textit{empty cluster} is a set of more than one third of the vertices of a complex such that no three of the vertices form a triangle and any proper superset of it contains a triangle.  We can specify that it contains $j$ vertices by calling it an \textit{empty cluster of size $j$}.}

\textup{An \textit{empty core} is an intersection of one or more empty clusters of the same size such that the intersection contains at least two vertices and any other empty cluster of the same size or larger contains at most one vertex of the empty core.}

\textup{An \textit{empty union} is a union of all empty clusters that contain a particular empty core and are of the same size as the empty clusters used to define the empty core.}

\textup{A \textit{column leg} is an intersection of two empty unions.  For each vertex of a column leg, one can count the number of triangles containing that vertex and no other vertex of the column leg.  The vertices of the column leg contained in tied for the most such triangles form a \textit{column foot}.}
\end{definition}

The terminology of the last paragraph may seem peculiar.  It is chosen to be descriptive in that a leg contains a foot, and both are a set of vertices at the bottom of a column, as shown in Lemmas~\ref{colleg} and \ref{colfoot}.

These were defined as they were because it is clear from the definitions that they are preserved by isomorphisms.  That is, an isomorphism between two complexes must send empty cores to empty cores, and so forth.  The preceding invariants give us ways to show that two complexes are not isomorphic, for example, if the sizes of their column legs do not match.

\begin{definition}
\textup{An \textit{indistinguishable pair} of vertices is a pair of distinct vertices $a$ and $b$ such that for all vertices $c$ and $d$ distinct from $a$ and $b$, the vertices $acd$ form a triangle exactly if $bcd$ do.  Equivalently, the vertex map swapping $a$ and $b$ induces an isomorphism of the complex.  An \textit{indistinguishable set} of vertices is a set of vertices that are pairwise indistinguishable.}
\end{definition}

We break up the proof that the complexes are not isomorphic into a series of lemmas.  We start by computing some small cases.

\begin{lemma} \label{basecase}
If $k \leq 2$, then no two constructions are isomorphic.
\end{lemma}

\proof  If $k \leq 1$, then there is only one construction with any particular number of vertices.  If $k = 2$ and $n = 3k+1 = 7$, there are three possible complexes.  One can count the number of vertices contained in exactly 11 triangles in each of these constructions and get totals of zero, one, and two.

If $n = 3k+2 = 8$, there are six possible complexes.  Among these six complexes, one has five empty clusters of size 4 (second row all red), one has four empty clusters of size 4 (second row all blue), and the rest have three empty clusters of size 4.  The last four complexes each have exactly two vertices contained in exactly one empty cluster of size 4.  If we count the number of triangles containing these two vertices for the various complexes, we get totals of 15 and 14, 14 and 14, 14 and 13, and 13 and 13.  \endproof

Lemma~\ref{basecase} leaves only the cases where $k \geq 3$.  As such, for the rest of this section, we assume that $k \geq 3$.  The next two lemmas show that most empty clusters are closely related to color sets.

\begin{lemma}  \label{emptyexist}
Any set of vertices contained in two columns such that all but the lowest in the left column are red and all but the lowest in the right column are blue has no triangles among any three of the vertices.
\end{lemma}

\proof  If three vertices are in the same column, they form a missing triangle of type 1.  If two vertices are in the left column and one is in the right, the higher of the two in the left is red, so they form a missing triangle of type 2.  If two vertices are in the right column and one is in the left, the higher of the two in the right column is blue, so they form a missing triangle of type 3.  \endproof

\begin{lemma} \label{emptysettypes}
Let $S$ be an empty cluster of vertices.  $S$ has vertices in exactly two columns.  All vertices of $S$ in a column except the lowest are the same color; we can refer to this as the color of the column of $S$.  If there are at least two vertices in each column, the left column is red and the right column is blue.

If the left column is blue or the right column red, the other column has exactly one vertex, it is the opposite color, it is the highest vertex of $S$, and $S$ has exactly $k+1$ vertices.  We call such a set a \textit{backwards empty cluster}.
\end{lemma}

\proof  If $S$ has vertices in all three columns, then $S$ contains a triangle consisting of one vertex from each column.  If all vertices of $S$ are in one column, we can add a vertex to $S$ in the column to the right if $S$ is red or to the left if $S$ is blue.  Hence, $S$ has vertices in exactly two columns.

If there are two vertices other than the lowest in the same column of $S$ that are different colors, then these two vertices and the lowest in the column form a triangle.

Suppose that both columns of $S$ have at least two vertices.  If both columns are red, then two vertices from the right column and the top one from the left form a triangle.  If both are blue, then two vertices from the left column and the top one from the right form a triangle.  If the left is blue and the right red, then the top vertex from each column and one other from the column with the highest vertex form a triangle.  This leaves only the possibility that the left column is red and the right blue.

If the left column is blue, then for its top two vertices not to form a triangle with the vertex from the right column, the vertex from the right column must be red and the highest of $S$.  The top row contains only red vertices, so the left column contributes at most $k$ vertices to $S$, and so $S$ has at most $k+1$ vertices.

Similarly, if the right column is red, then its top two vertices and the vertex from the left column form a triangle unless the vertex from the left column is blue and the highest of $S$.  The blue vertex cannot be in the top row, so the red column cannot have a vertex of $S$ in the top row.  Thus, $S$ has at most $k+1$ vertices in both cases.  $S$ must have at least $k+1$ vertices to be an empty cluster, so $S$ has exactly $k+1$ vertices. \endproof

Now we can show how to easily compute several of the invariants from the structure of Construction~\ref{subconst}.

\begin{lemma}  \label{k+2char}
There are no empty clusters of size $k+3$ or larger.  All empty clusters of size $k+2$ consist of a color set and one additional vertex in each column that is below the color set.  Conversely, any set of $k+2$ vertices that fits this description is an empty cluster.
\end{lemma}

\proof  By Lemma~\ref{emptysettypes}, for an empty cluster to have more than $k+1$ vertices, its left column must be red and its right column blue.  By Lemma~\ref{redbluesplit}, a color set has size either $k-1$ or $k$.  As such, the largest empty clusters possible are of size $k+2$ and are obtained by adding one additional vertex in each column, as in the statement of the lemma.

For the converse, by Lemma~\ref{emptyexist}, there are no triangles among the $k+2$ vertices.  As shown in the previous paragraph, the $k+2$ vertices cannot be a proper subset of a larger empty cluster.  \endproof

As an empty cluster must have at least $k+1$ vertices by definition, all empty clusters are of size $k+1$ or $k+2$.  By Lemma~\ref{redbluesplit}, each color set has size $k-1$ or $k$.  Hence, we can give them convenient labels.

\begin{definition}
\textup{A \textit{small empty cluster} is an empty cluster of size $k+1$.  A \textit{large empty cluster} is an empty cluster of size $k+2$.  A \textit{small empty core} is an empty core defined by an intersection of small empty clusters.  A \textit{large empty core} is an empty core defined by an intersection of large empty clusters.  A \textit{small color set} is a color set of size $k-1$.  A \textit{large color set} is a color set of size $k$.}
\end{definition}

\begin{lemma} \label{k+1char}
The small empty clusters are precisely the following constructions.
\begin{enumerate}
\item  Pick a small color set and add one additional vertex in each column that is below the color set.
\item  Pick a large color set, discard its lowest vertex in a column, add a vertex that is above the discarded vertex but below the next lowest vertex of the color set in that column, and add a vertex in the other column that is below the color set.
\item  Pick a backwards empty cluster.
\end{enumerate}
\end{lemma}

\proof  The third type is trivially an empty cluster.  The others have no triangles by Lemma~\ref{emptyexist}.  They cannot have vertices added to get a larger empty cluster by Lemma~\ref{k+2char}.

The third construction trivially includes all backwards empty clusters, so we only need to consider other empty clusters.  By Lemma~\ref{emptysettypes}, the empty cluster must involve exactly two columns, of which the left column must be red and the right column blue.  If the color set associated with the pair of columns is small, then all of these must be in the empty cluster in order to get $k+1$ vertices total.  Adding another vertex at the bottom of each column gives precisely the first construction.

Otherwise, the two columns have a large color set.  If we include all $k$ of these vertices in an empty cluster, then by Lemma~\ref{emptysettypes}, any additional vertices must be below the color set.  Adding one vertex does not prevent adding one in the other column, so this is not a small empty cluster.

By Lemma~\ref{emptysettypes}, the empty cluster can have at most two vertices outside of the color set.  As such, we can exclude only one of the $k$ vertices.  We then add two additional vertices below the lowest in their columns that remain to complete the prospective empty cluster.  The new vertex in the column with the excluded vertex must be above the excluded vertex to prevent adding it back, which gives precisely the second construction.  \endproof

\begin{lemma}  \label{emptycore}
The empty cores are precisely the sets consisting of a color set and the bottom vertex of any column for which the second lowest vertex is also in the color set, except for the cases where
\begin{enumerate}
\item $n = 3k+1$ and the unique blue vertex is the second row of the column with the highest vertex and higher than the red vertex to its right, or
\item $n = 3k+2$ and all vertices outside of the top row are blue.
\end{enumerate}
\end{lemma}

In subsequent lemmas, we refer to these last two constructions as the \textit{exceptional constructions}.

\proof  First suppose that the empty core is small and contained in a backwards empty cluster.  By Lemma~\ref{emptysettypes}, each empty cluster uses exactly two columns, so an empty core uses at most two columns.  Suppose that the empty core is contained in only one column.  To have more than one vertex, it must be the column with all but one vertex of the backwards empty cluster.

Suppose that the $k$ vertices in one column of a backwards empty cluster are all red.  The backwards empty cluster must have a blue vertex in the column to the left above all $k$ of these red vertices.  If there are two or more blue vertices, then at most one column has $k$ red vertices outside of the top row, and this column thus has a blue vertex to its right.  By Lemma~\ref{k+2char}, the $k$ red vertices are all contained in a large empty cluster, so the backwards empty cluster cannot be used to define an empty core.  The blue vertex cannot be in the top row, leaving only the possibility of exactly one blue vertex, which is the highest vertex of the second row.  This works in the case of $n = 3k+1$ but not $n = 3k+2$.

There is also the case where the $k$ vertices in a column of the backwards empty cluster are all blue.  If there is a red vertex in the column to its left, then by Lemma~\ref{k+2char}, the $k$ blue vertices are all contained in a large empty cluster.  All columns have a red vertex if $n = 3k+1$, and both columns with a vertex in the top row have a red vertex if $n = 3k+2$.  If there is a third red vertex in this latter case, it must be in the remaining column, leaving only the case where all vertices outside of the top row are blue.

Otherwise, the empty core has vertices from both columns.  Suppose that the empty core is contained in more than one backwards empty cluster.  If they do not have the same color of the colored column, then consecutive columns would be all red in one column and all blue in another, apart from the top row, which the construction does not allow.  If the backwards empty clusters have the same column with $k$ vertices, then the one vertex from the other column must vary, so the empty core only uses one column.

This leaves the possibility that the empty core is an intersection of a backwards empty cluster with a normal empty cluster in the same two columns.  Any vertex in the intersection must be in the backwards empty cluster, so it is blue if in the left column and red if in the right.  By Lemma~\ref{emptysettypes}, to be in the normal empty cluster also, it must be the lowest vertex in its column, and the only one of the opposite color. The only possible other higher vertices of the proper color which could be in the normal empty cluster are in the top row, which is all red, so at most one vertex there can be of the proper color.  As such, the normal empty cluster has at most three vertices.  Since $k \geq 3$, an empty cluster needs at least four vertices, a contradiction.

The other case is that the empty core does not have a backwards empty cluster used to define it.  By Lemma~\ref{emptysettypes}, each empty cluster uses exactly two columns.  If the empty clusters do not use the same two columns, then their intersection must be in a single column.  Suppose without loss of generality that the second lowest vertex in that column is red.  Then an empty cluster for which that column is blue contains one non-blue vertex, and an empty cluster for which the column is red contains no blue vertices.  The intersection contains at most one vertex, and does not contain an empty core.  Thus, all empty clusters containing the empty core use the same pair of columns.

Suppose that the empty core is large.  By Lemma~\ref{k+2char}, all large empty clusters include the entire color set for the pair of columns.  If a column has only one vertex below the color set, we are forced to use that vertex and it will be in the empty core.  If a column has more than one vertex below the color set, then different empty clusters will use different vertices, and their intersection contains none of the vertices at the bottom of the column.

Otherwise, the empty core is small.  If the color set is large, then the empty clusters are of type 2 as described in Lemma~\ref{k+1char}.  If any vertex neither in nor below the color set is in the empty core, then all vertices of the color set other than the one below this vertex and in its column must also be in the empty core.  This includes at least two vertices of the color set, so we can intersect with a large empty cluster and have at least two vertices.  If every vertex of the empty core is either in or below the color set, then a large empty cluster contains the entire empty core.  In either case, the empty core is not small.

The final possibility is that the color set is small.  This case works out the same as the one with large empty clusters:  all small empty clusters in these two columns must contain the entire color set, so its vertices are in the intersection.  The bottom vertex of each column is in the intersection exactly if it is in all small empty clusters of in the two columns, which happens precisely if the vertex immediately above it is of the column's color.  This is an intersection of all empty clusters in the two columns, so it does constitute an empty core.  \endproof

\begin{lemma}  \label{emptyunions}
If not an exceptional construction, an empty union consists of a color set and all vertices below the color set.
\end{lemma}

\proof  By Lemma~\ref{emptycore}, for each pair of columns, each color set is contained in an empty core.  By Lemma~\ref{emptysettypes}, we can fill these out to the largest empty clusters possible by adding one additional vertex at the bottom of each. We can choose any vertex below the last red one in the left column or the last blue one in the right column, so the union of all such empty clusters gives us the empty union as described in the lemma.  \endproof

\begin{lemma}  \label{colleg}
If not an exceptional construction, a set of $j \geq 2$ vertices forms a column leg exactly if they are all in the same column, the bottom $j$ vertices in the column, all the same color, and the $(j+1)$-th vertex from the bottom either does not exist or is the opposite color.
\end{lemma}

\proof  By Lemma~\ref{emptyunions}, there is exactly one empty union for each pair of columns.  An intersection of two empty unions thus has exactly one column in common.  The intersection consists of all vertices in the column which are either in or below both color sets.  This excludes the vertices that are above a vertex in the column of the opposite color, leaving the statement of the lemma.  \endproof

The next lemma is what allows us to disregard the exceptional constructions from the statement of some other lemmas.

\begin{lemma}  \label{exceptional}
No exceptional construction is isomorphic to any other construction.
\end{lemma}

\proof  Two exceptional constructions cannot be isomorphic to each other, as no two have the same number of vertices.

Suppose that $n = 3k+1$.  We can compute that the exceptional construction has column legs of sizes $k+1$, $k$, and $k-1$.  The only way for any other construction to have a column leg of size $k+1$ is for the entire column with the top vertex to be red, in which case, all vertices are red.  This construction has column legs of sizes $k+1$, $k$, and $k$, so it is not isomorphic to the exceptional one.

Now suppose that $n = 3k+2$.  We can compute that the exceptional construction has column legs of sizes $k+1$, $k$, and $k$.  In order for a normal construction to have a column leg of $k+1$ vertices, an entire column with one of the top two vertices must be red.  The other column with one of the top two vertices must then have a column leg of size $k$, so all but the top vertex must be blue.  As $k \geq 3$, this gives one column with at least two blue vertices and another with zero, which is impossible.  \endproof

\begin{lemma}  \label{colfoot}
If not an exceptional construction, a column foot consists of all vertices of a column leg that are indistinguishable from the lowest vertex of the column.
\end{lemma}

\proof  By Lemma~\ref{colleg}, all vertices of the column leg are the same color.  Suppose without loss of generality that they are all blue.  Let $u$ and $v$ be two vertices of the column leg that are distinguishable, and assume without loss of generality that $u$ is higher than $v$.  There must be two vertices $a$ and $b$ such that either $abu$ or $abv$ is a triangle, but not both.  Assume without loss of generality that $a$ is higher than $b$.

Let $C$ be the column containing $u$ and $v$, $D$ be the column to the right of $C$, and $E$ be the column to the left of $C$.  The following table lists the various possibilities for the columns of $a$ and $b$.  The entries of each row are the column of $a$, the column of $b$, and the conditions for $abu$ to form a triangle.

$$\begin{array}{ccl}
a & b & \textup{conditions} \\ C & C & a \textup{ and } b \textup{ different colors} \\ C & D & a \textup{ blue} \\ C & E & a \textup{ red} \\ D & C & (a \textup{ red and } a \textup{ below } u) \textup{ or } (a, b \textup{ both blue}) \\ D & D & a \textup{ red and } a \textup{ above } u \\ D & E & \textup{always} \\ E & C & a \textup{ and } b \textup{ both red} \\ E & D & \textup{always} \\ E & E & a \textup{ blue}
\end{array}$$

If $abu$ forming a triangle does not depend on $u$, then $abu$ and $abv$ either both form a triangle or else neither do.  We can thus restrict our attention to the $DC$ and $DD$ lines.  The only way for $abu$ to form a triangle and $abv$ not is if $a$ is red, in column $D$, below $u$, and above $v$, and $b$ is below $a$ and in column $C$, as given on line $DC$.  Note that this makes $b$ part of the column leg.  The only way for $abv$ to form a triangle and $abu$ not is if $a$ is red, in column $D$, below $u$, and above $v$, and $b$ is below $a$ and in column $D$, as given on line $DD$.

In both cases, in order for $u$ and $v$ to be distinguishable, there must be a red vertex below $u$, above $v$, and in the column $D$.  If there is such a red vertex below $u$, above $v$, and in column $D$, then we can take $a$ to be this vertex and $b$ to be the bottom vertex of column $D$, and get that $abv$ is a triangle while $abu$ is not, so $u$ and $v$ are distinguishable.  The bottom vertex of a column is below all colored vertices, so the set of vertices indistinguishable from the bottom one of a column leg is precisely the set of vertices of the column leg below all red vertices in the column to the right.

Note that only line $DD$ gives a triangle with only one vertex in the column leg.  For this case, vertices below all red vertices in the column to the right are in more triangles than vertices that are not.  By the definition of the column foot, the vertices it contains are those which are below all red vertices in the column to the right.  As we have seen, this is precisely the set of vertices indistinguishable from the lowest in the column.  \endproof

Now we have a way to essentially pick out the bottom row of a construction.  The basic idea of the proof is that we remove the bottom row and proceed by induction.

\begin{lemma} \label{row-1}
If two constructions are isomorphic, then if the bottom row of each were removed, they would still be isomorphic.
\end{lemma}

\proof  By Lemma~\ref{exceptional}, we only need to consider the constructions that are not among the exceptional cases.  It is clear from the definition of column feet that they are preserved by isomorphisms.  Pick a vertex from each column foot of one construction.  The isomorphism sends these to a vertex from each column foot of the other construction.  By Lemma~\ref{colfoot}, a vertex of a column foot is either the lowest vertex of the column or indistinguishable from it.  If the latter, there is an automorphism of the complex that interchanges the chosen vertex of the column foot with the lowest one of the column.  Composing these as necessary with the isomorphism between the two complexes, we get an isomorphism that preserves the bottom row.  This induces an isomorphism between the two complexes with their respective bottom rows removed.  \endproof

\begin{lemma} \label{emptyk+2}
If not an exceptional construction, there is an empty core containing $k+2$ vertices if and only if the second bottom row is of mixed color.  There is at most one empty core with $k+2$ vertices.
\end{lemma}

\proof  By Lemma~\ref{emptycore}, each empty core is a color set plus at most two vertices.  In order to contain $k+2$ vertices, it must be a large color set and have two additional vertices.  Each pair of columns has an associated empty core, and the bottom vertex of the column is in the empty core if and only if the second lowest vertex is.  This requires the second lowest vertex in the left column to be red and in the right column to be blue.  For this to happen more than once would require at least four vertices in the second lowest row.  Furthermore, it cannot happen if the second lowest row is solid color.

If the second bottom row has vertices of both colors, then at some point it has a red vertex with a blue vertex to its right.  The pair of columns of these two vertices has a corresponding color set.  By Lemma~\ref{redbluesplit}, if we remove the bottom two rows, at least $k-2$ vertices of the color set remain.  If we remove the bottom row, at most $k$ vertices of the color set remain.  Since the second bottom row adds two such vertices, both of these bounds are sharp, so the color set has $k$ vertices.  By Lemma~\ref{emptycore}, adding the two vertices in the bottom row gives an empty core with $k+2$ vertices.  \endproof

The above lemma is also true for the exceptional constructions, as backwards empty clusters do not affect large empty cores.

\begin{lemma}  \label{noniso}
No two constructions are isomorphic.
\end{lemma}

\proof  By Lemma~\ref{exceptional}, we can disregard the exceptional cases.  We use induction on $k$.  The base case of $k = 2$ is Lemma~\ref{basecase}.

For the inductive step, by Lemma~\ref{row-1}, if two complexes are isomorphic, then removing their bottom row of each leaves two isomorphic complexes.  By the inductive hypothesis, this requires the new complexes to be the same construction.  As such, only their bottom row of each can differ, so in the original constructions, only the second bottom rows can differ.  It thus suffices to show that if two complexes are identical except for their second bottom rows, then they are not isomorphic.

We can use Lemma~\ref{emptyk+2} to distinguish between the second bottom row being of solid color or mixed colors.  We can distinguish between it being all red or all blue by counting the number of column legs of size at least three; one color will have at least two such column legs, and the other at most one.

If the number of blue vertices in the second and third rows from the bottom added together is even, then the number of column legs of length at least three is odd; otherwise it is even.  This can distinguish between a complex with one blue vertex in the second bottom row and one with two blue vertices there.

That leaves as the only possible pairs of isomorphic constructions the cases where the vertex colorings are exactly the same, but the second bottom row has a blue vertex to the left of a red vertex, and the two complexes differ on which of the vertices is higher.  These two constructions differ only in a single triangle, which consists of the blue and red vertices already mentioned, plus the vertex in the bottom row of one column or the other.

The second bottom row is of mixed color, so by Lemma~\ref{emptyk+2}, there is an empty core that consists of a large empty cluster.  The two columns for this empty core are the same for both constructions in question, as by Lemma~\ref{emptycore}, empty cores depend only on the colors of the vertices.  They are not the same pair of columns as the triangle that varies is in, as the second bottom row has two vertices of this triangle, of which the left is blue and the right is red.  Hence, the empty core of size $k+2$ has one column in common with the triangle that varies between the two constructions.

The empty core of size $k+2$ contains the bottom vertex in both of its columns.  The triangle that varies between the two constructions has the bottom vertex in one of its two columns, and which column varies by construction.  Thus, the triangle has one vertex of the empty core of size $k+2$ in one construction and two vertices in the other.  The number of triangles with exactly two vertices in common with the empty core of size $k+2$ then differs by one between the two constructions, so they are not isomorphic.  \endproof

Finally, we can combine the above lemmas to prove Theorem~\ref{mytheorem}.

\smallskip\noindent {\it Proof of Theorem~\ref{mytheorem}: \ }  We have counted that Construction~\ref{subconst} gives ${1 \over 2}(6)^{k-1}$ complexes if $n = 3k+1$ and $6^{k-1}$ complexes if $n = 3k+2$.  By Theorem~\ref{optcolor}, all of them attain the bound of Tur\'{a}n's conjecture.  By Lemma~\ref{noniso}, no two of them are isomorphic.  \endproof

\section{Concluding remarks}

Kostochka observed that many constructions for the cases $n = 3k+1$ and $n = 3k+2$ could be obtained by removing one or two vertices from a construction for $n = 3(k+1)$.  Removing two vertices often gives a complex which does not attain the bound of Conjecture~\ref{turanconj}, however, and there can be many ways to remove one or two vertices and obtain the same complex.

Kostochka did not compute how many complexes he could construct in this manner.  From the results of this paper, one can readily show that he had $k2^{k-3}$ complexes if $n = 3k+2$ and $(k^2+3k-2)2^{k-4}$ complexes if $n = 3k+1$.  These are both far shy of the roughly $6^k$ complexes we have found.

Everything in Construction~\ref{subconst} is, however, a subcomplex of one given by Kostochka.  The easiest way to show this is from Construction~\ref{mainconst} by giving each vertex its own row, and filling in the empty spots in a row with a vertex of the same color as the vertex already in the row.  This can be done more efficiently to start with a smaller complex of Kostochka's construction and remove fewer vertices, but usually requires removing more than one fourth of the vertices.

Another question is whether these are all of the complexes that attain the bound of Tur\'{a}n's (3, 4)-conjecture.  If one could list all such complexes and show that none permit another triangle to be added, that would prove the conjecture.  Indeed, it is easy to show that adding one more triangle to any complex in Construction~\ref{subconst} would create a set of four vertices with four triangles.

The answer is that Construction~\ref{subconst} does not give all of the possible complexes that attain the bound, but it might be very close.  One can show by a brute force computer search that Construction~\ref{subconst} does give all of the complexes that attain the bound on eight or nine vertices, but it misses exactly one such complex if there are seven vertices.

The other complex is easy enough to describe:  use Construction~\ref{somecond} on six vertices, all colored red.  Add one additional vertex, with the triangles containing it precisely those that consist of the new vertex, one vertex from the top row, and one vertex from the bottom row.

This complex does not contain an empty cluster of size 4.  By the pigeonhole principle, any complex of \ref{mainconst} has a color set with two vertices that are not the bottom vertex of a column.  Adding the bottom vertex of each column gives an empty cluster on four vertices.  Note that any induced subcomplex of Construction~\ref{mainconst} is also a complex of this construction.  As such, not only is this new complex not isomorphic to any of Construction~\ref{subconst}, but it is not a subcomplex of any of them.

Neither is this new complex a subcomplex of some other complex that on at most twelve vertices attains the bound of the conjecture, as it is easy to show that every such complex must have a subcomplex on eight vertices that attains the bound, and we have the complete list of such complexes on eight vertices.  It could yet lead to a large class of additional complexes on much larger numbers of vertices, but it is quite plausible that this one new complex is simply an oddball exception.

Perhaps more likely is the possibility that there is some other large class of complexes which reduce to the known cases on nine or fewer vertices.  Note that Brown \cite{brown} had all of the constructions on nine vertices, even before Kostochka generalized his construction.

\end{document}